# On existence of general solution of the Navier-Stokes equations for 3D non-stationary incompressible flow

**Sergey V. Ershkov**

Institute for Time Nature Explorations,

M.V. Lomonosov's Moscow State University,

Leninskie gory, 1-12, Moscow 119991, Russia

e-mail: sergej-ershkov@yandex.ru

A new presentation of general solution of Navier-Stokes equations is considered here. We consider equations of motion for 3-dimensional non-stationary incompressible flow. The field of flow velocity as well as the equation of momentum should be split to the sum of two components: an irrotational (*curl-free*) one, and a solenoidal (*divergence-free*) one.

The obviously irrotational (*curl-free*) part of equation of momentum used for obtaining of the components of pressure gradient. As a term of such an equation, we used the irrotational (*curl-free*) vector field of flow velocity, which is given by the proper potential (besides, the continuity equation determines such a potential as a harmonic function).

The other part of equation of momentum could also be split to the sum of 2 equations: - with zero curl for the field of flow velocity (*viscous-free*), and the proper Eq. with viscous effects but variable curl. A solenoidal Eq. with viscous effects is represented by the proper Heat equation for each component of flow velocity with variable curl. Non-viscous case is presented by the PDE-system of 3 linear differential equations (in regard to the time-parameter), depending on the components of solution of the above Heat Eq. for the components of flow velocity with variable curl.

So, the existence of the general solution of Navier-Stokes equations is proved to be the question of existence of the proper solution for such a PDE-system of linear equations.

Final solution is proved to be the sum of 2 components: - an irrotational (*curl-free*) one and a solenoidal (*variable curl*) components.

1. **Introduction, the Navier-Stokes system of equations.**

In accordance with [1-3], the Navier-Stokes system of equations for incompressible flow of Newtonian fluids should be presented in the Cartesian coordinates as below (under the proper initial conditions):

$$\nabla \cdot \vec{u} = 0, \qquad (1.1)$$

$$\frac{\partial \vec{u}}{\partial t} + (\vec{u} \cdot \nabla)\vec{u} = -\frac{\nabla p}{\rho} + \nu \cdot \nabla^2 \vec{u} + \vec{F}, \qquad (1.2)$$

- where $u$ is the flow velocity, a vector field; $\rho$ is the fluid density, $p$ is the pressure, $\nu$ is the kinematic viscosity, and $F$ represents body forces (*per unit of mass in a volume*) acting on the fluid and $\nabla$ is the del (nabla) operator. Let us also choose the Ox axis coincides to the main direction of flow propagation.

Besides, we assume here external force $F$ above to be the *central* force, which has a potential $\phi$ represented by $F = -\nabla \phi$.

2. **The originating system of PDE for Navier-Stokes Eqs.**

Using the identity $(u \cdot \nabla)u = (1/2)\nabla(u^2) - u \times (\nabla \times u)$, we could present the Navier-Stokes equations in the case of incompressible flow of Newtonian fluids as below [4-5]:

$$\nabla \cdot \vec{u} = 0,$$

$$\frac{\partial \vec{u}}{\partial t} = \vec{u} \times \vec{w} + \nu \cdot \nabla^2 \vec{u} - \left(\frac{1}{2}\nabla(u^2) + \frac{\nabla p}{\rho} + \nabla \phi\right) \qquad (2.1)$$

- here we denote *the curl field* $w$, a pseudovector field [6] (*time-dependent*).

The system of equations (2.1) is equivalent to the Navier-Stokes system of equations for incompressible Newtonian fluids (1.1)-(1.2) in the sense of existence and smoothness of a general solution.

Let us denote as below (*according to the Helmholtz fundamental theorem of vector calculus*):

$$\nabla \times \vec{u} \equiv \vec{w}, \qquad \vec{u} \equiv \vec{u}_p + \vec{u}_w,$$

$$\nabla \cdot \vec{u}_w \equiv 0, \qquad \nabla \times (\vec{u}_p) \equiv 0,$$

- where $u_p$ is *an irrotational* (*curl-free*) field of flow velocity, and $u_w$ - is *a solenoidal* (*divergence-free*) field of flow velocity which generates a curl field $w$:

$$\vec{u}_p \equiv \nabla \varphi, \qquad \vec{u}_w \equiv \nabla \times \vec{A},$$

- here $\varphi$ - is the proper scalar potential, $A$ – is the appropriate vector potential. For such a potentials, we could obtain from the equation (1.1) the equality below

$$\nabla \cdot (\nabla \varphi + \nabla \times \vec{A}) = 0, \quad \Rightarrow \quad \Delta \varphi = 0, \tag{2.2}$$

- it means that $\varphi$ - is the proper *harmonic function* [6], time-dependent.

Thus, using *the curl of the curl* identity $\nabla \times (\nabla \times u) = \nabla(\nabla \cdot u) - \nabla^2 u$ and Eq. (2.2), the 2-nd equation of (2.1) could be presented as the sum of two equations below:

$$\begin{cases} \dfrac{1}{2} \nabla \{(\vec{u}_p + \vec{u}_w)^2\} + \dfrac{\nabla p}{\rho} + \nabla \phi = 0, \\ \\ \dfrac{\partial (\vec{u}_p + \vec{u}_w)}{\partial t} = (\vec{u}_p + \vec{u}_w) \times (\nabla \times \vec{u}_w) + \nu \cdot \nabla^2 \vec{u}_w, \end{cases} \tag{2.3}$$

- where the 1-st Eq. of (2.3) is the analogue of *Bernoulli* invariant [2].

So, if we solve the second equation of (2.3) in regard to the components of vectors ***u**** w*, ***u**** p* we could substitute it into the 1-st equation of (2.3); thus, we could obtain the proper expression for vector function $\nabla p$:

$$\frac{\nabla p}{\rho} = -\nabla \phi - \frac{1}{2}\nabla\{(\vec{u}_p + \vec{u}_w)^2\}, \qquad (2.4)$$

- where $\phi$ - is the also time-dependent function, in general case.

Thus, the key equation of Navier-Stokes equations in the case of incompressible flow for Newtonian fluids should be presented in a form below:

$$\frac{\partial(\vec{u}_p + \vec{u}_w)}{\partial t} = (\vec{u}_p + \vec{u}_w) \times (\nabla \times \vec{u}_w) + \nu \cdot \nabla^2 \vec{u}_w \qquad (2.5)$$

### 3. <u>Irrotational vs. solenoidal field of flow velocity in Navier-Stokes Eqs.</u>

Let us assume that some stable (constant) vorticity suddenly arise in *viscid* newtonian fluid; then, using *the curl of the curl* identity $\nabla^2 \boldsymbol{u} = \nabla(\nabla \cdot \boldsymbol{u}) - \nabla \times (\nabla \times \boldsymbol{u})$, we should obtain that $\nabla^2 \boldsymbol{u} = 0$ for such a flow (*viscous force equals to zero*).
It means that in the case of *viscid* newtonian fluids, the changing of *variable curl field **w*** strongly depends on *viscosity* factor of fluids; besides, the constant vorticity is zero.
But zero stable vorticity does not mean the zero velocity field which is assumed to be corresponding to such a curl field.

That's why Eq. (2.5) should be presented as the sum of two equations below:

$$\begin{cases} \dfrac{\partial \vec{u}_p}{\partial t} = \vec{u}_p \times \vec{w} + \vec{f}, & (3.1) \\[2ex] \dfrac{\partial \vec{u}_w}{\partial t} = \nu \cdot \nabla^2 \vec{u}_w, & (3.2) \end{cases}$$

- where we designate:

$$\vec{f} = \vec{u}_w \times \vec{w}.$$

Equation (3.2) immediately yields:

$$\frac{\partial \vec{w}}{\partial t} = \nu \cdot \nabla^2 \vec{w},$$

- it means indeed that the changing of *variable curl field* *w* should strongly depend on *viscosity* factor of fluids in case of *viscid* newtonian fluids.

So, if we solve Eq. (3.2) for the components of vector $\boldsymbol{u}_w$ (for finding of *w*), we should substitute it into Eq. (3.1) for solving it in regard to the vector function $\boldsymbol{u}_p$, which should depend on components of vector function $\boldsymbol{f}(\boldsymbol{u}_w, \boldsymbol{w})$ and curl *w*.

## 4. Field of flow velocity with zero curl.

Let us search for solutions $\boldsymbol{u}_p = \{U, V, W\}$ of (3.1), where each function is assumed to be depending on variables (x, y, z, *t*). It means that we should obtain time-dependent solutions to the system of PDE below:

$$\frac{\partial U}{\partial t} = (V \cdot w_z - W \cdot w_y) + f_x,$$

$$\frac{\partial V}{\partial t} = (W \cdot w_x - U \cdot w_z) + f_y, \qquad (4.1)$$

$$\frac{\partial W}{\partial t} = (U \cdot w_y - V \cdot w_x) + f_z,$$

- which should have a zero curl:

$$\frac{\partial W}{\partial y} - \frac{\partial V}{\partial z} = 0, \quad \frac{\partial U}{\partial z} - \frac{\partial W}{\partial x} = 0, \quad \frac{\partial V}{\partial x} - \frac{\partial U}{\partial y} = 0.$$

The "zero curl" condition yields the proper expressions for functions *U*, *V*, *W*:

$$U = \left(\frac{\partial \varphi}{\partial x}\right), \quad V = \left(\frac{\partial \varphi}{\partial y}\right), \quad W = \left(\frac{\partial \varphi}{\partial z}\right) \quad (4.2)$$

In such a case, the expressions (4.2) let us represent the system of Eqs. (4.1) as below:

$$\frac{\partial \left(\frac{\partial \varphi}{\partial x}\right)}{\partial t} = \left(\frac{\partial \varphi}{\partial y}\right) \cdot w_z - \left(\frac{\partial \varphi}{\partial z}\right) \cdot w_y + f_x ,$$

$$\frac{\partial \left(\frac{\partial \varphi}{\partial y}\right)}{\partial t} = \left(\frac{\partial \varphi}{\partial z}\right) \cdot w_x - \left(\frac{\partial \varphi}{\partial x}\right) \cdot w_z + f_y , \quad (4.3)$$

$$\frac{\partial \left(\frac{\partial \varphi}{\partial z}\right)}{\partial t} = \left(\frac{\partial \varphi}{\partial x}\right) \cdot w_y - \left(\frac{\partial \varphi}{\partial y}\right) \cdot w_x + f_z ,$$

- where

$$\left(\frac{\partial \varphi}{\partial x}\right), \quad \left(\frac{\partial \varphi}{\partial y}\right), \quad \left(\frac{\partial \varphi}{\partial z}\right),$$

- are the 3 unknown functions (time-dependent), which to be determined.

Thus, each of system of Eqs. (4.1) or (4.3) is proved to be equivalent to the system of Eqs. (3.1), but we should especially note that the combined system of equations (2.2), (2.4)+(3.1)+(3.2) is equivalent to the initial system of Navier-Stokes (1.1)+(1.2) in the sense of existence and smoothness of a solution.

## 5. Solution for the zero curl flow, $u_p = \{U, V, W\}$.

Eqs. (4.1) are the system of 3 linear PDEs (in regard to the time-parameter $t$) for 3 unknown functions: $U$, $V$, $W$, each of them could have the solution of *complex* value in general case [6]. Such a solution is out of direct physical significance but providing a new information about the equations themselves [7-8].

Besides, (4.1) is the system of 3 linear differential equations with all the coefficients depending on time $t$. In accordance with [6] p.71, the general solution of such a system should be given as below:

$$\chi_p(t) = \sum_{v=1}^{3} \zeta_{v,p} \cdot \left( \int \left( \frac{\Delta_v}{\Delta} \right) dt + C_v \right), \quad (p=1,2,3), \quad \Delta = \begin{vmatrix} \zeta_{1,1} & \zeta_{1,2} & \zeta_{1,3} \\ \zeta_{2,1} & \zeta_{2,2} & \zeta_{2,3} \\ \zeta_{3,1} & \zeta_{3,2} & \zeta_{3,3} \end{vmatrix} \quad (5.1)$$

- where $\{\chi_p\}$ - are the fundamental system of solutions of Eqs. (4.1): $\{U, V, W\}$, in regard to the time-parameter $t$; $\{\zeta_{v,p}\}$ - are the fundamental system of solutions of *the corresponding homogeneous* variant of (4.1), $\{C_v\}$ - are the set of functions, not depending on time $t$; besides, here we note as below:

$$\Delta_1 = \begin{vmatrix} f_x & f_y & f_z \\ \zeta_{2,1} & \zeta_{2,2} & \zeta_{2,3} \\ \zeta_{3,1} & \zeta_{3,2} & \zeta_{3,3} \end{vmatrix}, \quad \Delta_2 = \begin{vmatrix} \zeta_{1,1} & \zeta_{1,2} & \zeta_{1,3} \\ f_x & f_y & f_z \\ \zeta_{3,1} & \zeta_{3,2} & \zeta_{3,3} \end{vmatrix}, \quad \Delta_3 = \begin{vmatrix} \zeta_{1,1} & \zeta_{1,2} & \xi_{1,3} \\ \zeta_{2,1} & \zeta_{2,2} & \zeta_{2,3} \\ f_x & f_y & f_z \end{vmatrix}.$$

It means that the system of Eqs. (4.1) could be considered as having been solved if we obtain a general solution of *the corresponding homogeneous* system (4.1). Let us search for such a general solution of the corresponding homogeneous system (4.1) as below

$$\frac{\partial U}{\partial t} = V \cdot w_z - W \cdot w_y,$$

$$\frac{\partial V}{\partial t} = W \cdot w_x - U \cdot w_z, \qquad (5.2)$$

$$\frac{\partial W}{\partial t} = U \cdot w_y - V \cdot w_x.$$

The system of Eqs. (5.2) has *the analytical* way to present a proper solution [6] p.539. First of all, (5.2) yields the proper 1-st integral (invariant):

$$U^2 + V^2 + W^2 = \gamma^2(x, y, z)$$

- where $\gamma(x, y, z)$ – some arbitrary function, given by the initial conditions.

If we solve (5.2) for $\gamma = 1$, all other solutions could be obtained by the multiplying on the arbitrary function $\gamma(x, y, z)$. Thus, it is sufficient for us to consider the case $\gamma = 1$. If we consider a new functions $\xi(t)$, $\eta(t)$ for the case $\gamma = 1$ and, besides, we consider it under assumption below:

$$U + i \cdot V = \xi \cdot (\gamma - W), \quad U - i \cdot V = \eta^{-1} \cdot (W - \gamma), \qquad (*)$$

- where $i$ is the imaginary unit, then for $\xi(t)$ we should obtain the *Riccati* equation [6], [9]:

$$\xi' = \left(\frac{w_y + i \cdot w_x}{2}\right) \cdot \xi^2 - i \cdot w_z \cdot \xi + \left(\frac{w_y - i \cdot w_x}{2}\right) \qquad (5.3)$$

- but for $\eta(t)$ we should obtain the same as (5.3) *Riccati* equation (by the differentiating of the expression for a new function $\eta(t)$ above):

$$\eta' = \frac{(w_y + i \cdot w_x)}{2} \cdot \eta^2 - i \cdot w_z \cdot \eta + \frac{(w_y - i \cdot w_x)}{2}.$$

For the components $\{U, V, W\}$ the expressions above yield as below ($\xi \neq \eta \neq 0$):

$$U = \frac{(\gamma - W)\cdot(\xi - \eta^{-1})}{2}, \quad V = -\frac{(\gamma - W)\cdot i \cdot (\xi + \eta^{-1})}{2},$$

$$W = \gamma \cdot \frac{\left(1 + \dfrac{\eta}{\xi}\right)}{\left(1 - \dfrac{\eta}{\xi}\right)},$$

(5.4)

- but for the reason (*) we should note that:

$$\eta^{-1} = -\overline{\xi},$$

- that's why all the components $\{U, V, W\}$ (5.4) are *the real* functions in any case.

Besides, according to (1.1) the appropriate restriction should be valid for identifying of the function $\gamma(x, y, z)$ and the set of functions $\{C_v(x, y, z)\}$ in (5.1):

$$\frac{\partial U}{\partial x} + \frac{\partial V}{\partial y} + \frac{\partial W}{\partial z} = 0 \qquad (5.5)$$

- which is the PDE-equation of the 1-st kind; $w_x$, $w_y$, $w_z$ depend on variables (x,y,z, *t*).

So, the existence of the solution for Navier-Stokes system of equations (2.2), (2.4)+(3.1)+(3.2) (*which is equivalent to the initial system of Navier-Stokes* (1.1)+(1.2) *in the sense of existence and smoothness of a solution*) is proved to be the question of existence of the proper function $\gamma(x, y, z)$ and the set of functions $\{C_v(x, y, z)\}$ in (5.1) of so kind that the PDE-equation (5.5) should be satisfied under the given initial conditions (see [10], p.176). Besides, *non-homogeneous* solution (5.1) from (5.4) should generate a zero curl, according to the chosen form of solution (4.2); it means an additional restrictions at choosing of functions $\gamma(x, y, z)$ and $\{C_v(x, y, z)\}$ in (5.1):

$$\frac{\partial W}{\partial y} - \frac{\partial V}{\partial z} = 0, \quad \frac{\partial U}{\partial z} - \frac{\partial W}{\partial x} = 0, \quad \frac{\partial V}{\partial x} - \frac{\partial U}{\partial y} = 0.$$

## 6. Variable curl, the final presentation of solution.

Eq. (3.2) is known to be the Diffusion or Heat equation for each of components of flow velocity $u_w$ in the Cartesian coordinates [11]. Indeed, we know that the result of vector Laplacian for any vector function equals to the result of scalar Laplacian for each components of such a vector function only in the Cartesian coordinate systems [6].

Thus, we obtain that initial system of Navier-Stokes equations (1.1)-(1.2) for incompressible Newtonian fluids is equivalent to the combined system of equations (2.2), (2.4)+(3.1)+(3.2) in the sense of existence and smoothness of the solution.

To solve such a system, we present solution in a form $u = \nabla \varphi + u_w$, where *curl-free* field of flow velocity $\nabla \varphi = \{U, V, W\}$ is the flow velocity which generates a zero curl. Besides, $\varphi(x,y,z,t)$ is the proper *harmonic function* (potential) in regard to the variables $\{x,y,z\}$, see equation (2.2); but the flow velocity $u_w$ generates only a variable curl $w$ for the case of viscid newtonian fluid.

The 1-st, we should find the solution of equation (3.2) for each of components of flow velocity $u_w$; there are a modern methods for solving of 3D Heat equations [11]. The 2-nd, we should calculate the proper components of curl $w$. Then we should substitute all the components of curl $w$ and flow velocity $u_w$ to the equation (4.1) or (4.3) for solving it in regard to the components of flow velocity $\nabla \varphi = \{U, V, W\}$ in a form (4.2).

So, if we find a solution of flow velocity $u = \nabla \varphi + u_w$, we could substitute it into the equation (2.4) for the obtaining of a proper expression for vector function $\nabla p$.

Final solution is proved to be the sum of 2 components: - an irrotational (*curl-free*) one and a solenoidal (*variable curl*) components.

## 7. Discussions.

In fluid mechanics, a lot of authors have been executing their researches to obtain the analytical solutions of Navier-Stokes equations [12], even for 3D case of *compressible* gas flow [13].

We should especially note the excellent results [14]: it was outlined that the initial system of Navier-Stokes equations for incompressible 3D Newtonian fluids could be transformed to the Heat equation in regard to the proper components of the flow velocity, under the proper conditions (*if the initial conditions and the components of the externally applied force are chosen such that the nonlinear umbilical force vanishes by construction*). That's why such an equation has a proper analytical exact solution in R³ with the kernel of Gaussian type's integral [14].

Our presentation of the solution of Navier-Stokes equations for incompressible 3D fluids assumes implementing of the *Helmholtz* decomposition $\boldsymbol{u} = \nabla \varphi + \boldsymbol{u}_w$ (one of which generates *a zero* curl, other generates only *a variable* curl $\boldsymbol{w}$).

Besides, each component of the vector of flow velocity $\boldsymbol{u}_w$ (with *variable* curl) is proved evidently to be the solution of a proper Heat equation.

If we obtain a proper solutions for the flow velocity $\boldsymbol{u}_w$ (and curl $\boldsymbol{w}$) we could search for solutions of the flow velocity $\nabla \varphi = \{U, V, W\}$, with *zero* curl.

The system of equations for the flow velocity $\nabla \varphi = \{U, V, W\}$ is proved to have the proper solution of *real* value in general case. Having obtained the solution for flow velocity $\boldsymbol{u} = \nabla \varphi + \boldsymbol{u}_w$, it yields the proper expression for vector function of pressure gradient $\nabla p$.

We should especially note that the final solution is proved to be equivalent to the solution of the initial Navier-Stokes system of equations for incompressible Newtonian fluids in the sense of existence and smoothness of such a solution.

## 8. **Conclusion.**

A new presentation of general solution of Navier-Stokes equations is considered here. We consider equations of motion for 3-dimensional non-stationary incompressible flow. The field of flow velocity as well as the equation of momentum should be split to the sum of two components: an irrotational (*curl-free*) one, and a solenoidal (*divergence-free*) one.

The obviously irrotational (*curl-free*) part of equation of momentum used for obtaining of the components of pressure gradient. As a term of such an equation, we used the irrotational (*curl-free*) vector field of flow velocity, which is given by the proper potential (besides, the continuity equation determines such a potential as a harmonic function).

The other part of equation of momentum could also be split to the sum of 2 equations: - with zero curl for the field of flow velocity (*viscous-free*), and the proper Eq. with viscous effects but variable curl. A solenoidal Eq. with viscous effects is represented by the proper Heat equation for each component of flow velocity with variable curl. Non-viscous case is presented by the PDE-system of 3 linear differential equations (in regard to the time-parameter), depending on the components of solution of the above Heat Eq. for the components of flow velocity with variable curl.

So, the existence of the general solution of Navier-Stokes equations is proved to be the question of existence of the proper solution for such a PDE-system of linear equations. Final solution is proved to be the sum of 2 components: - an irrotational (*curl-free*) one and a solenoidal (*variable curl*) components.

**Acknowledgements**

I am thankful to CNews Russia project (Science & Technology Forum, mathematical branch) - for valuable discussions in preparing this manuscript. Especially I am thankful to Dr. L.Vladimirov-Paraligon for valuable discussions of this manuscript.